\documentclass[12pt]{article}
\usepackage{amsmath}
\usepackage{lineno}
\usepackage{txfonts}
\usepackage[utf8]{inputenc}
\usepackage[english]{babel}
\usepackage{hyperref}
\hypersetup{
	colorlinks=true,
	linkcolor=blue,
	filecolor=blue,
	urlcolor=blue,
	citecolor=blue,
}
\usepackage{booktabs}
\usepackage[pdftex]{graphics}
\graphicspath{{figures/}}
\usepackage{longtable}
\usepackage{amsfonts}

\usepackage{epstopdf}
\usepackage[section]{placeins}
\usepackage{adjustbox}
\usepackage[a4paper, total={6in, 8in}]{geometry}
\usepackage{cite}
\usepackage{breqn}
\usepackage{ltcaption}
\usepackage{framed}
\usepackage{caption}
\usepackage[font={scriptsize,it}]{caption}
\usepackage{epstopdf}
\usepackage{subcaption}
\usepackage{lipsum}
\usepackage{booktabs}
\usepackage{authblk}
\newcounter{mysection}
\usepackage{chngcntr}
\counterwithin{figure}{section}
\usepackage{titlesec}
\usepackage{hyperref}
\hypersetup{
	colorlinks=true,
	linkcolor=blue,
	filecolor=blue,
	urlcolor=blue,
	citecolor=blue,
}
\titleformat{\section}
{\normalfont\centering\large\scshape}{\thesection.}{1em}{}
\titleformat{\subsection}
{\it}{\thesubsection.}{1em}{}
\usepackage{graphicx}

\author{Awadh Pratap Singh$^a$ \footnote{\texttt{awadhma2015@gmail.com}}~, Shiv Prasad Yadav$^a$\footnote{ \texttt{spyorfma@gmail.com}}
	\\	{\small $^a$ Department of Mathematics, Indian Institute of Technology Roorkee, Roorkee-247667, India.}\\
}

\title{Performance evaluation of DMUs using hybrid fuzzy multi-objective data envelopment analysis}
\date{}

\begin{document}
	\maketitle
	\begin{abstract} 	The objective of this paper is to evaluate the performance of decision-making units (DMUs) using a hybrid fuzzy multi-objective (FMO) data envelopment analysis (DEA) approach. This study develops fuzzy multi-objective optimistic (FMOO) and pessimistic (FMOP) DEA models for performance evaluation of DMUs. To rank the DMUs, a ranking approach is used that can simultaneously integrate optimistic and pessimistic efficiencies. A comparison of the proposed approach with the existing approach is made with the help of an example. Finally, a real-life application of developed fuzzy optimistic and pessimistic DEA models in the education sector is presented.
	\end{abstract}
	
	\refstepcounter{mysection}
	\section{Introduction} \label{itm:intro7}
	Data envelopment analysis (DEA) is a non-parametric, data-aligned approach for assessing the performance of Decision-making Units (DMUs) that handle many inputs and outputs. Educational institutions, hospitals, banks, airlines, and other governmental agencies are examples of DMUs. Charnes, Cooper, and Rhodes \cite{charnes1978measuring} are known for inventing the DEA technique. The output-to-input ratio of a DMU is defined as efficiency (efficiency = output/input). The ratio of a DMU's efficiency to the largest efficiency under consideration is called relative efficiency. Its in the range of [0,1]. If a DMU's best relative efficiency equals to 1, it is considered efficient or optimistically efficient; otherwise, it is called non-efficient or optimistically non-efficient.\\
	
	DEA models assess efficiencies from both an optimistic and a pessimistic perspective. The worst relative efficiency, also known as pessimistic efficiency, is the efficiency computed from a pessimistic viewpoint. Efficiency is calculated as larger than 1 in the pessimistic DEA model. The DMUs are pessimistic in-efficient or non-inefficient if they have pessimistic efficiency values 1 or more than 1, respectively. We must assess both optimistic and pessimistic efficiencies simultaneously for the total performance of DMUs because we have two sorts of efficiencies: optimistic and pessimistic. In literature, Entani et al. \cite{entani2002dual}, Azizi H. \cite{azizi2011interval,azizi2014dea}  evaluated the efficiencies from both optimistic and pessimistic view points in a crisp environment while  Arya and Yadav \cite{arya2019development}  in the fuzzy environment. Gupta et al. \cite{gupta2020intuitionistic}  developed intuitionistic fuzzy optimistic and pessimistic multi-period portfolio optimization models. In their investigation, Puri and Yadav ]\cite{puri2015intuitionistic} created intuitionistic fuzzy optimistic, and pessimistic DEA models. The goal of these experiments was to create an interval using optimistic and pessimistic efficiency. Optimistic efficiency is the lower end of the interval in all of these studies, whereas pessimistic efficiency is the upper end.\\
	
	In all of the research mentioned above, the lower bound of optimistic efficiency and the upper bound of pessimistic efficiency were evaluated. The ranking is based on an interval created by combining the lower bound and upper bound efficiency of the optimistic and pessimistic efficiency DEA models. We are not in favor of considering only one bound for each optimistic and pessimistic efficiency. To overcome this shortcoming, we suggest considering both bounds of optimistic and pessimistic efficiency intervals for performance assessment. Based on this idea, we propose a fuzzy multi-objective optimistic (FMOO) and fuzzy multi-objective pessimistic DEA models to rank the DMUs. Awadh et al. \cite{awadhmulti2022} proposed a fuzzy multi-objective DEA model used to evaluate the performance of DMUs in a fuzzy environment. The advantage of the methodology is that it provides a uniform ranking of DMUs. There are several multi-objective optimization techniques exist in literature for solving DEA models \cite{wang2014fuzzy,chen2020comprehensive,boubaker2020role,zamani2020position}. But to the best of our knowledge, this is the first study using the fuzzy multi-objective technique for performance evaluation of DMUs in the optimistic and pessimistic environment. In this study, we developed FMOO and FMOP DEA models and used Wang et al's \cite{wang2007measuring} geometric average efficiency approach to rank the DMUs.\\
	
	The rest of the paper is organized as follows: Section 2 presents the preliminaries and some basic definitions. The proposed FMOO and FMOP DEA models along with the solving techniques are described in Section 3. The complete hybrid fuzzy multi-objective (FMO) DEA process is explained in Section 4. Section 5 includes the numerical illustration of the proposed methodology. Finally, in Section 6, the study's concluding remarks and future scope is presented.

	\section{Preliminary} \label{itm:prelim7}
	
	This section has given some important definitions and fuzzy operations that will help in developing Fuzzy optimistic and pessimistic DEA models.
	\begin{description}
		\item[Definition 1: \label{itm:Def 1}] [Fuzzy Set]\cite{zimmermann2011fuzzy} 	Let $\Omega$ be a universal set. A fuzzy set (FS) $\tilde{P}$ can be defined by
		
		$\tilde{P}=\{(\omega, \mu_{\tilde{P}}(\omega)):\omega\in \Omega\},\, where\,\,  \mu_{\tilde{P}}:\Omega\rightarrow{[0,1]}$.\\
		
		\item[Definition 2: \label{itm:Def 2}]  [Convex Fuzzy Set]\cite{zimmermann2011fuzzy}  A convex fuzzy set  $\tilde{P}$ is defined as if for all $\omega_{1}$, $\omega_{2}$ in	$\Omega$, 	\\
		$ min~\{\mu_{\tilde{P}}(\omega_{1}), \mu_{\tilde{P}}(\omega_{2})\} \leq \mu_{\tilde{P}}(\lambda \omega_{1}+(1-\lambda) \omega_{2}), \,\,\, where\,\, \lambda \in[0,1] .$\\
		
		\item[Definition 3: \label{itm:Def 3}] [Fuzzy number]\cite{zimmermann2011fuzzy} A fuzzy number (FN) $\tilde{P}$ is defined as a CFS $\tilde{P}$ of the real line $\mathbb{R}$ such that
		\begin{itemize}
			\item[$\mathrm{(i)}$] $\exists$ an unique $\omega_{o}\in{\mathbb{R}} \, with\, \mu_{\tilde{P}}(\omega_{o})=1;$
			
			\item[$\mathrm{(ii)}$] $\mu_{\tilde{P}}$ is piecewise continuous function.		
		\end{itemize}
		$\omega_o$ is called the modal value or mean value of $\tilde{P}.$\\

		\item[Definition 4: \label{itm: Def 4}] [Triangular Fuzzy Number]\cite{zimmermann2011fuzzy}
		The TFN $\tilde{P}$, denoted by $\tilde{P}=(p^L,p^M,p^U),$ (see Fig. \ref{fig:2.1}) is defined by the membership function $\mu_{\tilde{P}}$ given by
		\begin{equation}
		\mu_{\tilde{P}}(\omega) =
		\begin{cases}
		\frac{\omega-p^L}{p^M-p^L}, &\qquad p^L <\omega  \leq p^M;\\
		\frac{p^U-\omega}{p^U-p^M}, &\qquad p^M \leq \omega<p^U;\\
		0 & \text{otherwise.} \notag\\
		\end{cases}\\     
		\end{equation}
		$ \forall \omega \in \mathbb{R.}$ \\
		
		\begin{figure}[ht]
			\centering
			\includegraphics[width=0.7\linewidth]{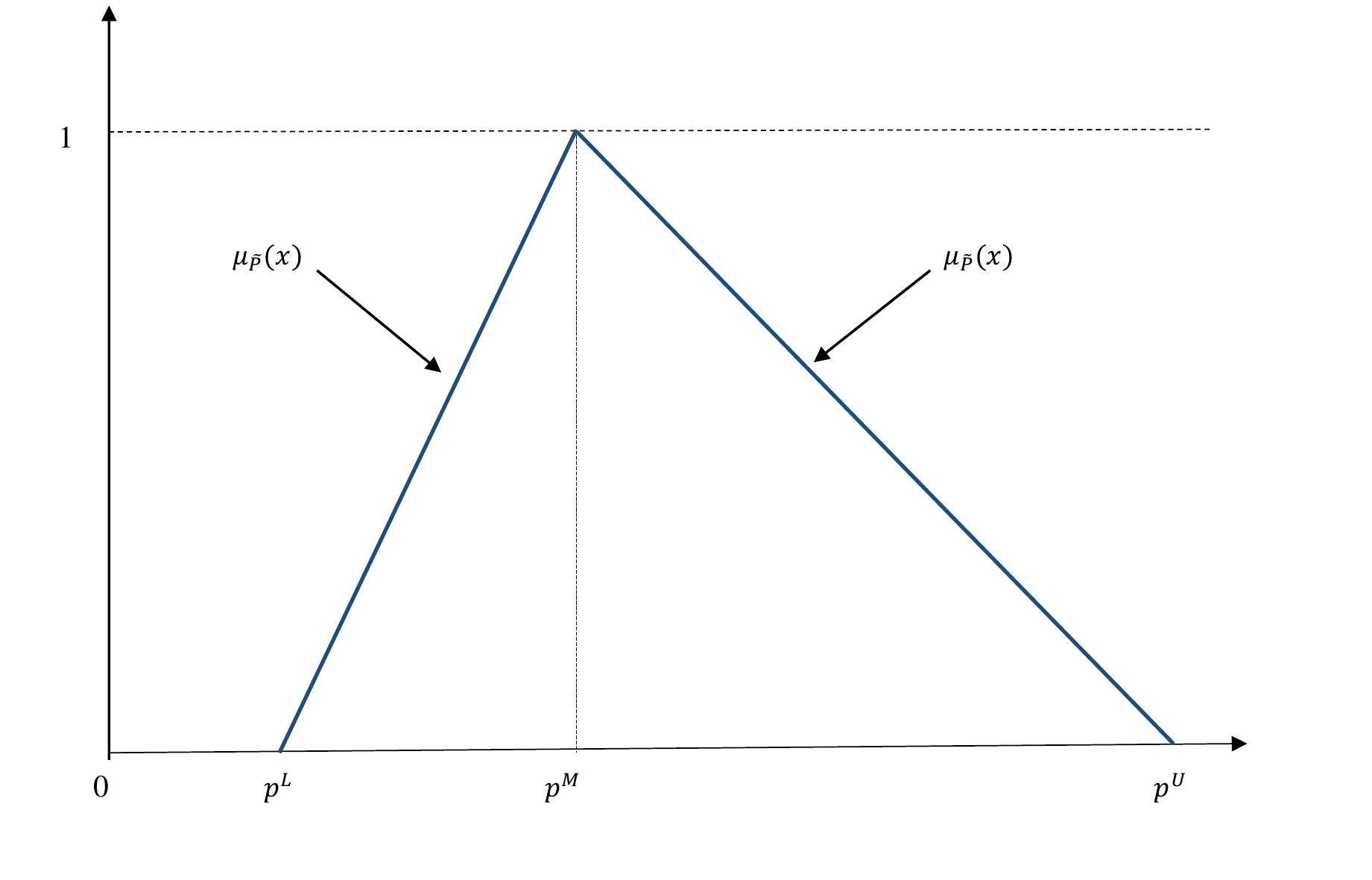}
			\caption{Triangular fuzzy number}
			\label{fig:2.1}
		\end{figure}
		
		\item[Definition 5: \label{itm: Def 5}] [Non-negative TFN]\cite{zimmermann2011fuzzy} A TFN $\tilde{P}=(p^L ,p^M, p^U)$ is said to be non-negative if and only if $p^L\geq0$.\\
		
		\item[Definition 6: \label{itm: Def 6}] [Positive TFN]\cite{zimmermann2011fuzzy} A TFN $\tilde{P}=(p^L ,p^M, p^U)$ is said to be positive if and only if $p^L>0$.\\
		
		\item[Definition 7: \label{itm: Def 7}] [$\alpha-cut$ ]\cite{zimmermann2011fuzzy} The $\alpha-cut$  of an FS $\tilde{P}$ in $\Omega$ is denoted by ${P}_{\alpha}$ and is defined by
		\begin{equation}
		{P}_{\alpha}=\{\omega\in \Omega: \mu_{\tilde{P}}(\omega)\geq\alpha\},\quad 0 \leq \alpha \leq 1.
		\end{equation}
		
		\textbf{Remark:} $P_{0}= \Omega.$
		
	\end{description}

	\subsection{\textbf{Arithmetic operations on TFNs}}
	Let $\tilde{P}=(p^{L},p^{M},p^{U})$ and $\tilde{Q}=(q^{L},q^{M},q^{U})$ be two positive TFNs. Then the aritmetic operations on TFNs are defined as follows\cite{wang2009fuzzy}:\\
	\begin{itemize}	
		\item[$\mathrm{(i)}$] \boldmath{Addition:} $\tilde{P}+\tilde{Q}=(p^{L}+q^{L},p^{M}+q^{M},p^{U}+q^{U} ),$
		\item[$\mathrm{(ii)}$] \boldmath{Subtraction:} $\tilde{P}-\tilde{Q}=(p^{L}-q^{U},p^{M}-q^{M},p^{U}-q^{L} ),$
		\item[$\mathrm{(iii)}$] \boldmath{Multiplication:} $\tilde{P}\times\tilde{Q} \approx (p^{L} q^{L},p^{M} q^{M},p^{U}  q^{U} ),$
		\item[$\mathrm{(iv)}$] \boldmath{Division:} $\tilde{P}/\tilde{Q} \approx (p^{L}/ q^{U},p^{M}/ q^{M},p^{U}/ q^{L} )$	
	\end{itemize}

	\section{Proposed Fuzzy Multi-objective Optimistic and Pessimistic DEA Models} \label{itm:propo7}
	Let us say, we wish to test the efficiency of n homogenous DMUs ($DMU_{j};j=1,2,3,...,n$). Assume that $DMU_j$ requires m inputs $ x_{ij},\,i=1,2,3,...,m$ to produce s outputs  $ y_{rj}, \,r=1,2,3,...,s$. Let $u_{ik}$ and $v_{rk}$ be the weights associated with $i^{th}$ input $x_{ik}$ and $r^{th}$ output $y_{rk}$ of $DMU_k\, (k=1,2,3,...,n).$ Let $E_k^{O}$ and $E_k^{P}$ stand for the optimistic and pessimistic efficiencies , respectively.  Entani et al. \cite{entani2002dual} proposed the optimistic and pessimistic DEA models given in Table  \ref{tab:Table-1}.
	
	\begin{table}[ht]
		\centering
		\caption{{\textbf{Optimistic and pessimistic DEA models}}}
		\label{tab:Table-1}
		\begin{tabular}{cc}
			\hline
			
			Optimistic DEA model {(Model 1)}\label{itm: model1}  & Pessimistic DEA model {(Model 2)}\label{itm: model2}\\\hline
			For $k=1,2,3,...,n,$ & For $k=1,2,3,...,n,$\\
			$\max E_{k}^O = \dfrac{\sum\limits_{r=1}^{s}y_{rk}{v_{rk}}}{\sum\limits_{i=1}^{m}x_{ik}{u_{ik}}}$ & $\min E_{k}^P = \dfrac{\sum\limits_{r=1}^{s}y_{rk}{v_{rk}}}{\sum\limits_{i=1}^{m}x_{ik}{u_{ik}}}$ \\
			{\hbox {subject to }}$ \dfrac{\sum\limits_{r=1}^{s}y_{rj}{v_{rk}}}{\sum\limits_{i=1}^{m}x_{ij}{u_{ik}}}\leq 1\,\, \forall j$, & {\hbox {subject to }} $ \dfrac{\sum\limits_{r=1}^{s}{y_{rj}v_{rk}}}{\sum\limits_{i=1}^{m}{x_{ij}u_{ik}}}\geq 1\,\, \forall j,$\\
			$u_{ik}\geq \varepsilon\,\,\forall i,\,\, v_{rk}\geq \varepsilon\,\,\forall r,\,\,\,\varepsilon > 0.$ & $u_{ik}\geq \varepsilon\,\,\forall i, \,\,   v_{rk}\geq \varepsilon\,\,\forall r,\,\,\forall k,\,\,\varepsilon > 0.$ \\
			\hline
		\end{tabular}
	\end{table}  
	
	\begin{description}
		\item[\textbf{Definition 8:}\label{itm: Def 8}\cite{arya2019development}] Let the optimal values of the optimistic and pessimistic DEA models for $DMU_k$ be $E_k^{O*}$ and $E_k^{P*}$, respectively. Then $DMU_k$ is said to be \textbf{optimistic efficient} if $E_k^{O*}=1;$ otherwise \textbf{optimistic non-efficient.} On the other hand, $DMU_k$ is said to be \textbf{pessimistic inefficient} if $E_k^{P*}=1;$ otherwise \textbf{pessimistic non-inefficient.} 
	\end{description}

	Due to the ambiguity and fluctuation of such real-world data, it is challenging to get accurate and reliable input and output data. Assume that the fuzzy inputs and outputs for the $DMU_j$,  $\tilde{x}_{ij}$ and $\tilde{y}_{rj}$, respectively. Then the fuzzy optimistic (FO) and fuzzy pessimistic (FP) DEA models are described (see Table \ref{tab:Table-2}) as follows: 
	
	\begin{table}[ht]
		\centering
		\caption{\textbf{FO and FP DEA models}}
		\label{tab:Table-2}
		\begin{tabular}{cc}
			\hline
			FO DEA model {(Model 3)}\label{itm: model3}& FP DEA model {(Model 4)}\label{itm: model4} \\\hline
			For $k=1,2,3,...,n,$ & For $k=1,2,3,...,n,$\\
			$\max \tilde{E}_{k}^O = \dfrac{\sum\limits_{r=1}^{s}\tilde{y}_{rk}{v_{rk}}}{\sum\limits_{i=1}^{m}\tilde{x}_{ik}{u_{ik}}}$ & $\min \tilde{E}_{k}^P = \dfrac{\sum\limits_{r=1}^{s}\tilde{y}_{rk}{v_{rk}}}{\sum\limits_{i=1}^{m}\tilde{x}_{ik}{u_{ik}}}$ \\
			{\hbox {subject to }}$ \dfrac{\sum\limits_{r=1}^{s}\tilde{y}_{rj}{v_{rk}}}{\sum\limits_{i=1}^{m}\tilde{x}_{ij}{u_{ik}}}\leq \tilde{1}\,\, \forall j$, & {\hbox {subject to }} $ \dfrac{\sum\limits_{r=1}^{s}{\tilde{y}_{rj}v_{rk}}}{\sum\limits_{i=1}^{m}{\tilde{x}_{ij}u_{ik}}}\geq \tilde{1}\,\, \forall j,$\\
			$u_{ik}\geq \varepsilon\,\,\forall i,\,\, v_{rk}\geq \varepsilon\,\,\forall r,\,\,\,\varepsilon > 0.$ & $u_{ik}\geq \varepsilon\,\,\forall i, \,\,   v_{rk}\geq \varepsilon\,\,\forall r,\,\,\forall k,\,\,\varepsilon > 0.$ \\
			\hline
		\end{tabular}
	\end{table}
	
	Assume that the fuzzy input $\tilde{x}_{ij}=(x_{ij}^{L}, x_{ij}^{M}, x_{ij}^{U})$, fuzzy output $\tilde{y}_{rj}=(y_{rj}^{L}, y_{rj}^{M}, y_{rj}^{U})$ for the $DMU_j$, and $\tilde{1}=(1,1,1)$ are taken as TFNs. Then FO and FP models can be transformed into triagular fuzzy optimistic (TFO) and triangular fuzzy pessimistic (TFP) DEA models as follows (see table \ref{tab:Table-3}):
	
	\begin{table}[ht]
		\centering
		\caption{\textbf{TFO and TFP DEA models}}
		\label{tab:Table-3}
		\begin{tabular}{cc}
			\hline
			TFO DEA model {(Model 5)}\label{itm: model5} & TFP DEA model {(Model 6)}\label{itm: model6}\\\hline
			For $k=1,2,3,...,n,$ & For $k=1,2,3,...,n,$\\
			$\max ({E}_{k}^{O,L}, {E}_{k}^{O,M}, {E}_{k}^{O,U}) = \dfrac{\sum\limits_{r=1}^{s}(y_{rk}^{L}, y_{rk}^{M},y_{rk}^{U}){v_{rk}}}{\sum\limits_{i=1}^{m}(x_{i}^{L}, x_{ik}^{M}, x_{ik}^{U}){u_{ik}}}$ & $\min ({E}_{k}^{P,L}, {E}_{k}^{P,M}, {E}_{k}^{P,U}) = \dfrac{\sum\limits_{r=1}^{s}(y_{rk}^{L}, y_{rk}^{M}, y_{rk}^{U}){v_{rk}}}{\sum\limits_{i=1}^{m}(x_{ik}^{L}, x_{ik}^{M}, x_{ik}^{U}){u_{ik}}}$ \\
			{\hbox {subject to }}$ \dfrac{\sum\limits_{r=1}^{s}(y_{rj}^{L}, y_{rj}^{M}, y_{rj}^{U}){v_{rk}}}{\sum\limits_{i=1}^{m}(x_{ij}^{L}, x_{ij}^{M}, x_{ij}^{U}){u_{ik}}}\leq (1,1,1)\,\, \forall j$, & {\hbox {subject to }} $ \dfrac{\sum\limits_{r=1}^{s}(y_{rj}^{L}, y_{rj}^{M}, y_{rj}^{U})v_{rk}}{\sum\limits_{i=1}^{m}(x_{ij}^{L}, x_{ij}^{M}, x_{ij}^{U}) u_{ik}}\geq (1,1,1)\,\, \forall j,$\\
			$u_{ik}\geq \varepsilon\,\,\forall i,\,\, v_{rk}\geq \varepsilon\,\,\forall r,\,\,\,\varepsilon > 0.$ & $u_{ik}\geq \varepsilon\,\,\forall i, \,\,   v_{rk}\geq \varepsilon\,\,\forall r,\,\,\forall k,\,\,\varepsilon > 0.$ \\
			\hline
		\end{tabular}
	\end{table}
	
	Now we will propose a methodology to solve TFO and TFP DEA models given in Table \ref{tab:Table-3}. In this methodology, we will use Awadh et. al's FMODEA \cite{awadhmulti2022}. First, let us try to develop FMOO DEA model for the efficiency evaluation of DMUs in an optimistic sense. The TFO DEA model, described in Table \ref{tab:Table-3} is re-written in Model 7.

	\begin{description}
		\item[\textbf{Model 7}\label{itm:model7}]  For $ k=1,2,3,...,n,$
		\begin{flalign}
		& Max \, (E^{O,L}_{k}, E^{O,M}_{k}, E^{O,U}_{k})= \frac{(\sum_{r=1}^{s}v_{rk}y^{L}_{rk},\sum_{r=1}^{s}v_{rk}y^{M}_{rk},\sum_{r=1}^{s}v_{rk}y^{U}_{rk})}{(\sum_{i=1}^{m}u_{ik}x^{L}_{ik},\sum_{i=1}^{m}u_{ik}x^{M}_{ik},\sum_{i=1}^{m}u_{ik}x^{U}_{ik})} &\\
		& {\hbox {subject to }}\frac{(\sum_{r=1}^{s}v_{rk}y^{L}_{rj},\sum_{r=1}^{s}v_{rk}y^{M}_{rj},\sum_{r=1}^{s}v_{rk}y^{U}_{rj})}{(\sum_{i=1}^{m}u_{ik}x^{L}_{ij},\sum_{i=1}^{m}u_{ik}x^{M}_{ij},\sum_{i=1}^{m}u_{ik}x^{U}_{ij})} \leq(1,1,1)\qquad\forall j=1,2,3,...,n; &\\ 
		&u_{ik},v_{rk}\geq \varepsilon\quad \forall i,r,j.&\notag
		\end{flalign}
	\end{description}
	
	Model 7 can be transformed into the Model 8 by using division rule of two fuzzy numbers as follows:
	
	\begin{description}
		\item[\textbf{Model 8}\label{itm:model8}]  For $ k=1,2,3,...,n,$
		\begin{flalign}
		& Max \, (E^{O,L}_{k}, E^{O,M}_{k}, E^{O,U}_{k})= \Big(\frac{\sum_{r=1}^{s}v_{rk}y^{L}_{rk}}{\sum_{i=1}^{m}u_{ik}x^{U}_{ik}},\frac{\sum_{r=1}^{s}v_{rk}y^{M}_{rk}}{\sum_{i=1}^{m}u_{ik}x^{M}_{ik}},\frac{\sum_{r=1}^{s}v_{rk}y^{U}_{rk}}{\sum_{i=1}^{m}u_{ik}x^{L}_{ik}}\Big) &\\
		& {\hbox {subject to }}\Big(\frac{\sum_{r=1}^{s}v_{rk}y^{L}_{rj}}{\sum_{i=1}^{m}u_{ik}x^{U}_{ij}},\frac{\sum_{r=1}^{s}v_{rk}y^{M}_{rj}}{\sum_{i=1}^{m}u_{ik}x^{M}_{ij}},\frac{\sum_{r=1}^{s}v_{rk}y^{U}_{rj}}{\sum_{i=1}^{m}u_{ik}x^{L}_{ij}}\Big) \leq(1,1,1)~~\forall j=1,2,3,...,n; \label{itm: eq.4}\\
		&u_{ik},v_{rk}\geq \varepsilon\quad \forall i,r,j.&\notag
		\end{flalign}
	\end{description}
	
	As TFNs have the following property (see Eq. \ref{itm: Eq.5})
	\begin{equation} \label{itm: Eq.5}
	\frac{\sum_{r=1}^{s}v_{rk}y^{L}_{rj}}{\sum_{i=1}^{m}u_{ik}x^{U}_{ij}} \leq \frac{\sum_{r=1}^{s}v_{rk}y^{M}_{rj}}{\sum_{i=1}^{m}u_{ik}x^{M}_{ij}} \leq \frac{\sum_{r=1}^{s}v_{rk}y^{U}_{rj}}{\sum_{i=1}^{m}u_{ik}x^{L}_{ij}} 
	\end{equation}
	
	Now using the property explained in eq. \ref{itm: Eq.5}, Model 8 can be transfomed into FMOO DEA model as follows:
	
	\begin{description}
		\item[\textbf{Model 9}\label{itm:model9}] \textbf{(Proposed FMOO DEA model):}  For $ k=1,2,3,...,n,$
		\begin{flalign}
		& Max~\Big[\frac{\sum_{r=1}^{s}v_{rk}y^{L}_{rk}}{\sum_{i=1}^{m}u_{ik}x^{U}_{ik}},\,\, \frac{\sum_{r=1}^{s}v_{rk}y^{M}_{rk}}{\sum_{i=1}^{m}u_{ik}x^{M}_{ik}},\,\, \frac{\sum_{r=1}^{s}v_{rk}y^{U}_{rk}}{\sum_{i=1}^{m}u_{ik}x^{L}_{ik}}\Big] &\\
		& {\hbox {subject to }}\frac{\sum_{r=1}^{s}v_{rk}y^{U}_{rj}}{\sum_{i=1}^{m}u_{ik}x^{L}_{ij}} \leq 1 \quad \forall j=1,2,3,...,n; \label{itm: eq.7}\\
		&u_{ik},v_{rk}\geq \varepsilon\quad \forall i,r,j.&\notag
		\end{flalign}
	\end{description}
	
	Now using Charnes-Cooper transformation \cite{charnes1997data}, Model 9 can be convereted into Model 10 as follows:
	
	\begin{description}
		\item[\textbf{Model 10}\label{itm:model10}] \textbf{(Proposed FMOP DEA model):}  For $ k=1,2,3,...,n,$
		\begin{flalign}
		& Min~\Big[{\sum_{r=1}^{s}v_{rk}y^{L}_{rk}},\,\, {\sum_{r=1}^{s}v_{rk}y^{M}_{rk}},\,\, {\sum_{r=1}^{s}v_{rk}y^{U}_{rk}}\Big] &\\
		& {\hbox {subject to }}{\sum_{i=1}^{m}u_{ik}x^{U}_{ik}}=1&\\
		&{\sum_{i=1}^{m}u_{ik}x^{M}_{ik}}=1&\\
		&{\sum_{i=1}^{m}u_{ik}x^{L}_{ik}}=1&\\			
		&	{\sum_{r=1}^{s}v_{rk}y^{U}_{rj}}-{\sum_{i=1}^{m}u_{ik}x^{L}_{ij}} \leq 0 \quad \forall j=1,2,3,...,n; \label{itm: eq.9}\\
		& u_{ik},v_{rk}\geq \varepsilon\quad \forall i,r,j.&\notag
		\end{flalign}
	\end{description}
	
	Hence, Model 10 is deterministic FMOO DEA model. Similarly, we will convert the TFP DEA model into FMOP DEA model by using Awadh et.al's \cite{awadhmulti2022} FMODEA model as follows:
	
	\begin{description}
		\item[\textbf{Model 11}\label{itm:model11}] \textbf{(Proposed FMOP DEA model):}  For $ k=1,2,3,...,n,$
		\begin{flalign}
		& Min~\Big[\frac{\sum_{r=1}^{s}v_{rk}y^{L}_{rk}}{\sum_{i=1}^{m}u_{ik}x^{U}_{ik}},\,\, \frac{\sum_{r=1}^{s}v_{rk}y^{M}_{rk}}{\sum_{i=1}^{m}u_{ik}x^{M}_{ik}},\,\, \frac{\sum_{r=1}^{s}v_{rk}y^{U}_{rk}}{\sum_{i=1}^{m}u_{ik}x^{L}_{ik}}\Big] &\\
		& {\hbox {subject to }}\frac{\sum_{r=1}^{s}v_{rk}y^{L}_{rj}}{\sum_{i=1}^{m}u_{ik}x^{U}_{ij}} \geq 1 \quad \forall j=1,2,3,...,n; \label{itm: eq.9}\\
		& u_{ik},v_{rk}\geq \varepsilon\quad \forall i,r,j.&\notag
		\end{flalign}
	\end{description}
	
	Now using Charnes-Cooper transformation \cite{charnes1997data}, Model 11 can be convereted into Model 12 as follows:
	
	\begin{description}
		\item[\textbf{Model 12}\label{itm:model12}] \textbf{(Proposed FMOP DEA model):}  For $ k=1,2,3,...,n,$
		\begin{flalign}
		& Min~\Big[{\sum_{r=1}^{s}v_{rk}y^{L}_{rk}},\,\, {\sum_{r=1}^{s}v_{rk}y^{M}_{rk}},\,\, {\sum_{r=1}^{s}v_{rk}y^{U}_{rk}}\Big] &\\
		& {\hbox {subject to }}{\sum_{i=1}^{m}u_{ik}x^{U}_{ik}}=1&\\
		&{\sum_{i=1}^{m}u_{ik}x^{M}_{ik}}=1&\\
		&{\sum_{i=1}^{m}u_{ik}x^{L}_{ik}}=1&\\			
		&	{\sum_{r=1}^{s}v_{rk}y^{L}_{rj}}-{\sum_{i=1}^{m}u_{ik}x^{U}_{ij}} \geq 0 \quad \forall j=1,2,3,...,n; \label{itm: eq.9}\\
		& u_{ik},v_{rk}\geq \varepsilon\quad \forall i,r,j.&\notag
		\end{flalign}
	\end{description}

	Hence, Model 12 is deterministic FMOP DEA model which provides the the efficiency of DMUs in pessimistic sense.

	\subsection{\textbf{Algorithm for solving proposed FMOO and FMOP DEA models}}\label{Sec. 4.2}
	The following is a summary of the step-by-step algorithm for solving the reduced MOO and FMOP DEA Model (i.e, Model 10 and 12) as follows:\\
	\begin{enumerate}
		
		\item[Step 1:]  For each DMU convert multi-objective problem into a single objective problem using weighted sum method \cite{marler2010weighted}.
		\item[Step 2:] Generate a population of random weights. Suppose that problem has $D$ objectives to optimize and $p$ variables. Then generate ($100\times p$) set of $D$ weights. There is no thumb rule for population selection. It depends upon the decision-maker. For computational purpose, we take hundred times of the number of variables present in the problem.

		\item[Step 3:] Solve the single objective problem for each $DMU$ generated in Step 1 with the help of weights generated in Step 2.

		\item[Step 4:] From step 3, we get ($100 \times p$) Pareto solutions. Find the best one among all solutions obtained.
		
		\item[Step 5:] Best weighted solution obtained from Step 4 can be treated as efficiency of DMU.
		
		\item[Step 6:] Use above steps to determine $E_{k}^{O*}$ and $E_{k}^{P*}$.

	\end{enumerate}

	After getting both optimistic and pessimistic efficiency scores we are needed to rank the DMUs by considering both the efficiencies simultaneously. To rank the DMUs we will use the geometric average efficiency approach, proposed by Wang et al. \cite{wang2007measuring}.  According to Wang et al. \cite{wang2007measuring}, if $E_{k}^{O*}$ and $E_{k}^{P*}$ are the optimistic and pessimistic efficiencies, respectively, for $DMU_k,$ then the geometric average efficiency ($E_{k}^{geometric}$) can be defined as follows:
	\begin{equation}
	E_{k}^{geometric}= \sqrt{E_{k}^{*O} \times E_{k}^{*P}} \label{Eq: 20}
	\end{equation}
	
	Wang et al. \cite{wang2007measuring} proposed this approach for the DMUs with crisp input and output data. We extend the idea for the DMUs with fuzzy input and output data, particularly triangular fuzzy data.		
	
	\section{\textbf{Complete hybrid FMO DEA performance decision process}}		
	
	In complex real world problems, the hybridization of DEA/FDEA/FMO DEA using other techniques is very effective \cite{entani2002dual,azizi2011interval,arya2019development}. These researchers handled both the optimistic and pessimistic DEA models simultaneously. The hybridization process can be divided into four steps as follows: 	  
	
	\begin{itemize}
		\item [(i)]\textbf{Input-output data selection and collection phase:} In this phase, the decision-maker follows the following steps: (a) selection of the relevant input and output data variables, (ii) Collection of input-output data quantitatively and qualitatively, (iii) classification of the data according to its nature, crisp or fuzzy, (iv) fuzzification of the data based on criteria and expert's opinion.\\
		
		\item[(ii)] \textbf{Efficiency measurement phase:} During this phase, the decision-maker chooses the best method for evaluating performance. The best and worst performance are obtained using the suggested fuzzy optimistic and pessimistic DEA models. As a result, the suggested FMO DEA strategy leads to the overall performance of DMUs. As a result, the hybrid FMO DEA technique is ideal for a fair decision-making process.\\ 
		
		\item[(iii)] \textbf{Ranking phase:} Using the geometric average ranking approach, the complete ranking is obtained by selecting optimistic and pessimistic FMO DEA models.\\
		
		\item[(iv)] \textbf{Recommendation phase:}This is the final stage of the decision-making process, in which policy-makers and experts are given recommendations based on the ranking results obtained in the previous phase. The recommendations include suggestions for critical modifications that the management must undertake in order to increase the DMUs' efficiencies.
	\end{itemize}
	
\begin{figure}
		\centering
		\includegraphics[width=0.9\linewidth]{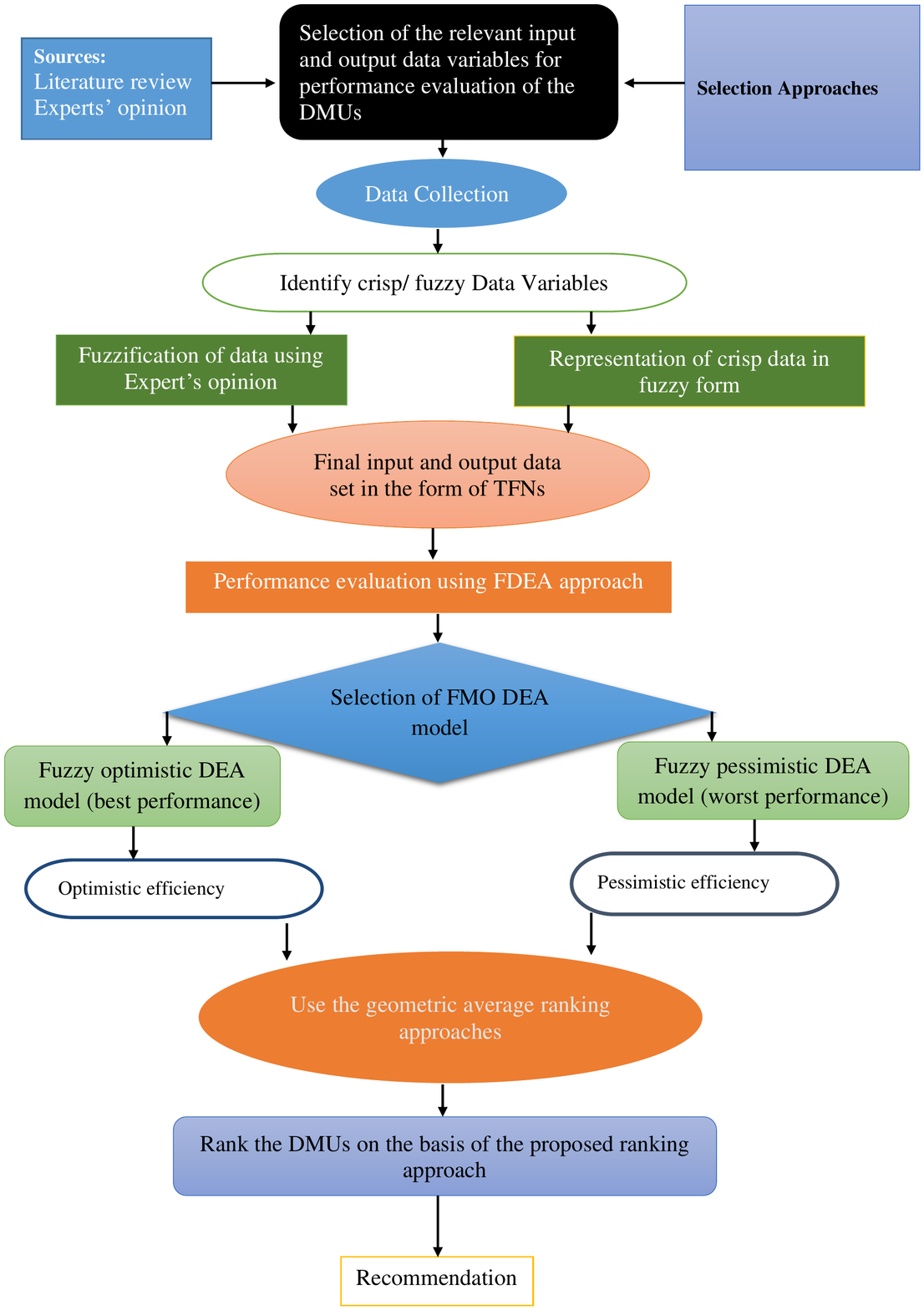}
		\caption{Schematic view of the proposed hybrid FMO DEA performance decision model}
		\label{fig: 1}
	\end{figure}

	Figure \ref{fig: 1}, depicts the suggested hybrid FDEA performance efficiency evaluation process. 
	
	\section{Numerical Illustration} 
	In this part, we use an example to demonstrate the efficacy of the presented models. This example is based on Guo and Tanaka's study \cite{guo2001fuzzy}, which used five DMUs with two fuzzy inputs and two fuzzy outputs. We try to validate the provided models and ranking approach using this scenario. A case study of the proposed models in the field of education is also provided.
	
	\subsection{\textbf{Numerical Illustration: An example}} 
	Table \ref{tab:Table 4} presents the fuzzy input-output data for the problem considered by Guo and Tanaka, which has five DMUs with two fuzzy inputs and two fuzzy outputs.  
	
	\begin{table}[ht]
		\centering
		\caption{{\textbf{Fuzzy input-output data for 5 DMUs} (Source: Guo and Tanaka \cite{guo2001fuzzy})}}
		\label{tab:Table 4}
		\begin{tabular}{ccccc}
			\hline
			
			DMUs & \multicolumn{2}{c}{Fuzzy inputs} & \multicolumn{2}{c}{Fuzzy outputs}\\ \hline
			\,\, & Input 1 & Input 2 & Output 1 & Output 2 \\ \hline
			$DMU_1$ & (3.5, 4.0, 4.5) &(1.9, 2.1, 2.3)& (2.4, 2.6, 2.8)&(3.8, 4.1, 4.4)\\
			$DMU_2$&  (2.9, 2.9, 2.9) & (1.4, 1.5, 1.6)&(2.2, 2.2, 2.2)& (3.3, 3.5, 3.7)\\
			$DMU_3$&(4.4, 4.9, 5.4)&(2.2, 2.6, 3.0)&(2.7, 3.2, 3.7)&(4.3, 5.1, 5.9)\\
			$DMU_4$&(3.4, 4.1, 4.8)&(2.1, 2.3, 2.5)&(2.5, 2.9, 3.3)&(5.5, 5.7, 5.9)\\
			$DMU_5$&(5.9, 6.5, 7.1)&(3.6, 4.1, 4.6)&(4.4, 5.1, 5.8)& (6.5, 7.4, 8.3)\\\hline
		\end{tabular}
	\end{table}
	
	The optimistic and pessimistic efficiency scores are calculated by using the proposed FMOO and FMOP DEA models respectively. Table \ref{tab:Table 5} presents the outcomes ($E_{k}^{*O}$ and $E_{k}^{*P}$) of both the recommended models. 
	
	\begin{table}
		\centering
		\caption{Optimistic efficiency scores ($E_{k}^{*O}$) of DMUs (Example) (Source: Guo and Tanaka \cite{guo2001fuzzy})}
		\label{tab:Table 5}
		\resizebox{\textwidth}{!}{
			\begin{tabular}{|l|lll|llll|l|}				
				\hline
				DMUs & \multicolumn{3}{l|}{ \quad Most favourable weights} & \multicolumn{4}{l|}{\qquad \qquad  Most favourable solutions} & Efficiency ($E_{k}^{*O}$)  \\ \hline
				& $~~w_1$      &$~~ w_2 $    &~~ $ w_3 $     &~~ $v_1$   & ~~ $v_2$   & $ ~~~~ u_1$   & ~~$u_2$   &         \\ \hline	
				$D_A$   &	0.0283	&0.0414&	0.9303&	0.1017 &	0.0858	&2.05E-05&	0.4733 &0.6579
				\\
				$D_B$   &	0.0065&	0.0003&	0.9932&	3.68E-05&	0.1987	&3.45E-01&	1.00E-05 &0.7347					\\			
				$D_C$  &	0.0183&	0.0315	&0.9502&	1.13E-01&	0.0459	&2.03E-01&	1.00E-05 &0.6822	
				\\	
				$D_D$  &	0.0094& 	0.1102&	0.8804&	0.0020&	1.36E-01&	1.00E-05&	4.33E-01 & 0.8061	
				\\				
				$D_E$   &	0.0238&	0.0189&	0.9573&	1.33E-01&	0.0040&	1.53E-01&	1.00E-05 & 0.8011
				\\			
				\hline			
			\end{tabular}
		}
	\end{table}

	\begin{table}
		\centering
		\caption{Pessimistic efficiency scores ($E_{k}^{*P}$) of DMUs (Example) (Source: Guo and Tanaka \cite{guo2001fuzzy})}
		\label{tab:Table 6}
		\resizebox{\textwidth}{!}{
			\begin{tabular}{|l|lll|llll|l|}				
				\hline
				DMUs & \multicolumn{3}{l|}{ \quad Most favourable weights} & \multicolumn{4}{l|}{\qquad \qquad  Most favourable solutions} & Efficiency ($E_{k}^{*P}$)  \\ \hline
				& $~~w_1$      &$~~ w_2 $    &~~ $ w_3 $     &~~ $v_1$   & ~~ $v_2$   & $ ~~~~ u_1$   & ~~$u_2$   &         \\ \hline	
				$D_A$   &	0.9702&	0.0136	&0.0162&	0.1056&	0.2639&	1.00E-05&	0.4733&	1.2611
				
				\\
				$D_B$   &	0.9711	&0.0066	&0.0223&	5.87E-05&	0.4330&	0.3448&	1.00E-05 &	1.4335
				\\			
				$D_C$  &	0.9257&	0.0444	&0.0298&	0.1905&	1.35E-01&	2.03E-01&	1.00E-05&1.1158
				
				\\	
				$D_D$  &	0.9462&	0.0401&	0.0138&	0.4805	&1.43E-04&	1.21E-05&	4.33E-01&1.2151
				
				\\				
				$D_E$   &		0.8932	&0.1068&	7.15E-05&	0.0004	&1.92E-01&	1.53E-01&	1.00E-05&1.2675	
				\\			
				\hline			
			\end{tabular}
		}
	\end{table}
	
	Based on the efficiency scores ($E_{k}^{*O}$, $E_{k}^{*P}$) obtained from FMOO, and FMOP DEA; the geometric efficiency score is obtained by using Eq (\ref{Eq: 20}). The ranking is done based on the $E_{k}^{geometric}$. The ranking obtained from the proposed FMO DEA model (see Table 7) is compared with Wang et al.'s \cite{wang2006centroids} method. Since both the methods provide different level of efficiency scores; so it will be appropriate to compare the ranks of DMUs with the help of Spearman's rank correlation coefficient \cite{zar2005spearman}. The correlation coefficient ($\rho=0.883$) indicates that efficiency obtained from both the models are highly correlated. This indicates that the proposed methodology is quite efficient to rank the DMUs in fuzzy environment.

	\begin{table}
		\centering
		\caption{  \textbf{The proposed and geometric average efficiency scores and ranks of DMUs} }
		\label{tab:Table 7}
			\begin{tabular}{|l|l|l|ll|ll|}				
				\hline
				DMUs &  $E_{k}^{*O}$& $E_{k}^{*P}$ & $E_{k}^{geometric}$ & Proposed rank & Wang et al.& Rank\\ \hline
				
				$D_A$   &0.6579&	1.2611 & 0.9108&~~~~~~~4&1.090&3
				\\
				$D_B$   &0.7347&1.4335&1.0262&~~~~~~~1&1.154&1
				\\			
				$D_C$  &0.6822&1.1158&0.8725&~~~~~~~5&1.092&2
				
				\\	
				$D_D$  &0.8061&1.2151&0.9867&~~~~~~~3&1.154&1	
				
				\\				
				$D_E$   &0.8011&1.2675&1.0076 &~~~~~~~2&1.154&1
				\\			
				\hline			
			\end{tabular}
		
	\end{table}
	
	\subsection{\textbf{Education sector application}} \label{itm:Sec 3.1}
	We investigate the application of checking the performance efficiency of Indian Institutes of Management (IIMs) in India to apply the given approach for checking efficiency. The Indian Institutes of Management (IIMs) are administrative and research-oriented educational institutes. They mostly provide undergraduate, graduate, doctoral, and other courses. The Ministry of Education of India has acknowledged IIMs to be the most important institutions in the country. This is a real-world application in which we used the following two inputs and two outputs for 13 IIMs:
	\begin{itemize}
		\item [(i)] \textbf{Input 1:} Number of students ($x_{1}$)
		\item [(ii)] \textbf{Input 2:} Number of faculty members ($x_{2}$)
		\item [(iii)] \textbf{Output 1:} Placements and higher studies ($y_{1}$)
		\item [(iv)] \textbf{Output 2:} Publications  ($y_{2}$)
	\end{itemize}
	
	The DMUs are the IIMs. The information was gathered from the NIRF website, which was introduced on September 29, 2015 by the Honorable Minister of Human Resource Development (The Ministry of Education). The data was collected during a four-year period, from 2015-16 to 2018-19. Then the collected data is having four entries for each input and output so it is converted into TFN representing only one entry as $\tilde{a}=(a^L, a^M, a^U)$ by following process.\\	
	$a^L=$ Minimum value of the data for any particular IIM\\
	$a^M=$ Average value of the data for any particular IIM\\
	$a^U=$ Maximum value of the data for any particular IIM\\
	
	Two fuzzy inputs are used: the number of students and the number of faculty members. In this study, the number of students placed or who moved on to further studies, as well as the number of publications, are considered two fuzzy outputs. On the basis of data collected from the NIRF website for four years, the data is transformed into TFNs. TFNs are used as input-output data in Table \ref{tab:Table-8}. The $E_{k}^{*O}$ and $E_{k}^{*P}$ is calculated by using the proposed FMOO DEA (Model 10) and FMOP DEA (Model 12) models respectively. The efficiencies are calculated by using the proposed algorithm (Sect. 3.1). The results for optimistic and pessimistic efficiencies are shown in Table 9 and Table 10 respectively. The decision-makers have complete control over the input and output data. The ranking will be altered if any decision-maker selects different input-output data sets. The goal of this research is to provide a novel approach for evaluating DMU performance. The ranking is done based on the efficiencies obtained from Table \ref{tab:Table 9} and Table \ref{tab:Table 10}. The proposed ranking method is used and the ranks of DMUs are presented in Table \ref{tab:Table 11}.
	
	\begin{table}
		\centering
		
		\caption{\textbf{Input and output data for IIMs}}
		\label{tab:Table-8}
		\resizebox{\textwidth}{!}{
			\begin{tabular}{|lllllll|}
				\hline
				DMU & IIM Name      & State   \qquad\qquad  & \multicolumn{2}{l}{Inputs} & \multicolumn{2}{l|}{Outputs} \\ \hline
				&               &           & $\qquad\tilde{x_1}$   & $\qquad\tilde{x_2}$  & $\qquad\tilde{y_1}$  &$\qquad \tilde{y_2}$  \\ \hline
				$D_1$  & IIM Bangalore & Karnataka&	(424,682,955)&	(91,104,113)  & (393,410,449) &	(72,136,212)        \\
				
				$D_2$  & IIM Ahmedabad & Gujrat &	(461,715,992)&	(91, 112,128) & (411,421,427)&	(30,109,217)         \\
				
				$D_3$  & IIM Calcutta & West Bengal &	(487,803,1042)&	(86,94,105)  &  (483,505,535)&	(29,109,207)       \\
				
				$D_4$  & IIM Lucknow & Uttar Pradesh &	(455,725,990)&	(81,88,95) &   (440,456,506)&	(8,65,126)       \\
				
				$D_5$  & IIM Indore & Madhya Pradesh &	(549,1020,1657)&	(73,94,104)  &  (508,593,634)&	(16,66,141)       \\
				
				$D_6$  & IIM Kozhikode & Kerala &	(370,593,806)&	(58,69,77) &     (347,360,382)&	(28,74,97)     \\
				
				$D_7$  & IIM Udaipur & Rajasthan &	(120,260,419)&	(21,46,101) &  (120,136,171)&	(14,37,67)        \\
				
				$D_8$  & IIM Tiruchirapalli & Tamilnadu &	(108,228,387)&	(25,37,52)  &  (102,121,172)&	(5,16,24)       \\
				
				$D_9$  & IIM Raipur & Chhatisgarh  &	(160,270,438) &	(21,33,46)& (111,140,193) &	(12,35,52)           \\
				
				$D_{10}$  & IIM Rohtak & Haryana & (158,276,428) &	(20,28,34) &   (137,146,155)&	(31,52,69)	       \\
				
				$D_{11}$  & IIM Shillong & Meghalaya  &	(155,263,365) &	(27,27,28)  &  (118,146,172)&	(3,15,30)       \\
				
				$D_{12}$  & IIM Kashipur & Uttarakhand &	(125,262,472)&	(15,28,38) &   (101,123,164)&	(7,21,36)       \\
				
				$D_{13}$  & IIM Ranchi & Jharkhand &	(189,300,452)&	(16,29,40)  &   (156,169,178)&	(11,22,51)      \\

				\hline
			\end{tabular}
		}
	\end{table}

	\begin{table}
		\centering
		\caption{\textbf{Optimistic efficiency scores ($E_{k}^{*O}$) of IIMs}}
		\label{tab:Table 9}
		\resizebox{\textwidth}{!}{
			\begin{tabular}{|l|lll|llll|l|}				
				\hline
				DMUs & \multicolumn{3}{l|}{ \quad Most favourable weights} & \multicolumn{4}{l|}{\qquad \qquad  Most favourable solutions} & Efficiency ($E_{k}^{*O}$)  \\ \hline
				& $~~w_1$      &$~~ w_2 $    &~~ $ w_3 $     &~~ $v_1$   & ~~ $v_2$   & $ ~~~~ u_1$   & ~~$u_2$   &         \\ \hline	
				$D1$   &	0.01787&	0.0778&	0.9043&	444.153&	3324.761&	2285.578&	2.33E+03 &0.7226 			
				\\
				$D2$   &	0.0478&	0.1005&	0.8518&	1396.474&	16176.274&	9418.478&	10935.333&0.6898
				\\			
				$D3$  &	0.0499&	0.196&	0.754&	806.660&	2433.738&	1982.503&	5583.592&	0.5587 				
				\\	
				$D4$  &	0.0221&	0.3803& 0.5976&	1546.952&	2440.449&	2880.0246&	9311.586038&0.4327	
				\\				
				$D5$   &	0.0477&	0.1005&	0.8518&	365.137&	820.407&	480.326&	9054.454&0.4108
				\\			
				$D6$	&0.2034&	0.0314&	0.7652&	1602.965&	3615.045&	3725.355&	20785.561&0.3215\\
				
				$D7$		&0.3551&	0.1704&	0.4744&	642.319&	1230.289&	1288.628&	2208.538 &0.5495\\
				$D8$	&	0.3664&	0.2204&	0.4132&	425565.259&	203161.842&	647291.776&	757192.279&0.4879\\
				$D9$		&0.0383&	0.1884&	0.7732&	338.6579&	1003.258&	701.589&	7803.547&0.3668\\
				$D10$	&	0.1080&	0.1712&	0.7207&	115.1186&	423.1365&	290.529&	3449.675&0.3456\\
				$D11$	&	0.2775&	0.0378&	0.6846&	375.977&	475.993&	644.032&	2112.405&0.3981\\
				$D12$	&	0.0650&	0.4640&	0.4709&	1410.191&	3379.327&	2175.496&	32052.563 &0.3053\\
				$D13$	&	0.0845&	0.3342&	0.5813&	283.888&	1052.874&	229.074&	10865.027&0.3498\\ 			
				\hline			
			\end{tabular}
		}
	\end{table}

	\begin{table}
		\centering
		\caption{\textbf{Pessimistic efficiency scores ($E_{k}^{*P}$) of IIMs}}
		\label{tab:Table 10}
		\resizebox{\textwidth}{!}{
			\begin{tabular}{|l|lll|llll|l|}				
				\hline
				DMUs & \multicolumn{3}{l|}{ \quad Most favourable weights} & \multicolumn{4}{l|}{\qquad \qquad  Most favourable solutions} & Efficiency ($E_{k}^{*O}$)  \\ \hline
				& $~~w_1$      &$~~ w_2 $    &~~ $ w_3 $     &~~ $v_1$   & ~~ $v_2$   & $ ~~~~ u_1$   & ~~$u_2$   &         \\ \hline	
				$D1$   &0.9356&	0.052&	0.0123&	0.0081&	1.01E-05&	1.00E-05&	0.0095&	3.2013\\
				$D2$&0.9683&	0.0217&	0.0099&	0.0074&	1.01E-05&	1.00E-05&	0.0088&	3.0711\\
				$D3$&	0.9848&	0.0006&	0.0144&	0.0087&	1.01E-05&	1.00E-05&	0.0103&	4.2401\\
				$D4$&	0.9077&	0.0701&	0.0221&	0.0004&	1.08E-01&	1.00E-05&	0.0112&	1.7888\\
				$D5$&	0.942&	0.0136&	0.0443&	0.00047&	1.03E-01&	1.00E-05&	0.0107&	2.5246\\
				$D6$ &0.9702&	0.0161&	0.0136&	0.0121&	1.01E-05&	1.00E-05&	0.0144&	4.2338\\
				$D7$&	0.9903&	0.0013&	0.0083&	0.0144&	3.55E-05&	3.10E-03&	1.00E-05&	1.7457\\
				$D8$&	0.971&	0.0065&	0.0223&	0.0065&	1.19E-01&	1.00E-05&	0.0241&	1.3322\\
				$D9$&0.9367&	0.0353&	0.0278&	0.023&	1.36E-04&	1.00E-05&	0.0273&	2.6345\\
				$D10$&0.0247&	0.0294&	2.25E-05&	1.00E-05&	0.0349&	4.033&	4.2989&	4.5642
				\\
				$D11$&	0.9437&	0.0369&	0.0192&	0.0098&	1.79E-01&	1.00E-05&	0.0364&	1.8937\\
				$D12$&	0.9515&	0.0105&	0.0378&	0.0277&	1.47E-04&	1.00E-05&	0.0329&	2.8733\\
				$D13$&	0.9538&	0.0406&	0.0054&	0.0131&	1.99E-05&	2.52E-03&	0.0036&	2.0613\\ 			 			
				\hline			
			\end{tabular}
		}
	\end{table}

	\begin{table}
		\centering
		\caption{  \textbf{The proposed geometric average efficiency scores and ranks of DMUs} }
		\label{tab:Table 11}
		\
		\begin{tabular}{|l|l|l|l|l|}				
			\hline
			DMUs &  $E_{k}^{*O}$& $E_{k}^{*P}$ & $E_{k}^{geometric}$ & Proposed rank \\ \hline
			
			$D1$&0.7226& 3.2013&1.5209 &~~~~~~~~2\\
			$D2$& 0.6898&3.0711&1.4554&~~~~~~~~3\\
			$D3$&0.5587&4.2401&1.5392&~~~~~~~~1\\
			$D4$&0.4328&1.7888&0.8798&~~~~~~~~10	\\
			$D5$&0.4108&2.5246&1.1084&~~~~~~~~6\\
			$D6$ &0.3215&4.2338&1.1667&~~~~~~~~5\\
			$D7$&0.5495&1.7457&0.9795&~~~~~~~~8\\
			$D8$&0.4879&1.3323&0.8063&~~~~~~~~13\\	
			$D9$&0.3668&2.6345&0.9829&~~~~~~~~7\\
			$D10$&0.3456&4.0468&1.1826&~~~~~~~~4\\
			$D11$&0.3981&1.8937&0.8683&~~~~~~~~11\\
			$D12$&0.3050&2.8733&0.9366&~~~~~~~~9	\\
			$D13$&0.3498&2.0613&0.8491&~~~~~~~~12	\\ 
			
			\hline			
		\end{tabular}
		
	\end{table}
	
	\section{Conclusions}
	To solve a DEA model with imprecise parameters, a multi-objective optimization strategy is proposed in this paper. The data for individual DMUs is fuzzified to model the imprecision in the data points and to aggregate numerous points into a single data input. The model is then translated into FMOO and FMOP DEA models employing arithmetic operations between TFNs, which are then used to quantify the performance efficiency of DMUs. DMUs are Indian Institutes of Management (IIMs), which are ranked according to their efficiency scores. Among the 13 IIMs included for this study, IIM Calcutta is the best performing institute, while IIM Tiruchirappalli is the worst performing. This model has the advantage of providing a uniform ranking of DMUs with fuzzy inputs and fuzzy outputs. This paper's key contribution is the development of a theoretical model of the DEA framework in the presence of an ambiguous data. The current study initiates the interest in  the fuzzy DEA, which is identified by uncertainty in the data points.\\
	
	Other DEA models that have been extended in unfavourable contexts may be studied in the future. By solving a MOOP model with different random weights, the current research aims to improve the neutrality of DEA models in testing related performance. However, the proposed method is confined to a few of the LPP solution hypotheses. This can also be examined in several types of uncertainty, such as stochastic and robust.
	
	\bibliographystyle{ieeetr}
	\bibliography{reference7}

\begin{thebibliography}{10}

\bibitem{charnes1978measuring}
A.~Charnes, W.~W. Cooper, and E.~Rhodes, ``Measuring the efficiency of decision
  making units,'' {\em European journal of operational research}, vol.~2,
  no.~6, pp.~429--444, 1978.

\bibitem{entani2002dual}
T.~Entani, Y.~Maeda, and H.~Tanaka, ``Dual models of interval dea and its
  extension to interval data,'' {\em European Journal of Operational Research},
  vol.~136, no.~1, pp.~32--45, 2002.

\bibitem{azizi2011interval}
H.~Azizi, ``The interval efficiency based on the optimistic and pessimistic
  points of view,'' {\em Applied Mathematical Modelling}, vol.~35, no.~5,
  pp.~2384--2393, 2011.

\bibitem{azizi2014dea}
H.~Azizi, ``Dea efficiency analysis: A dea approach with double frontiers,''
  {\em International Journal of Systems Science}, vol.~45, no.~11,
  pp.~2289--2300, 2014.

\bibitem{arya2019development}
A.~Arya and S.~P. Yadav, ``Development of optimistic and pessimistic models
  using fdea to measure performance efficiencies of dmus,'' {\em International
  Journal of Process Management and Benchmarking}, vol.~9, no.~3, pp.~300--322,
  2019.

\bibitem{gupta2020intuitionistic}
P.~Gupta, M.~K. Mehlawat, S.~Yadav, and A.~Kumar, ``Intuitionistic fuzzy
  optimistic and pessimistic multi-period portfolio optimization models,'' {\em
  Soft Computing}, pp.~1--26, 2020.

\bibitem{puri2015intuitionistic}
J.~Puri and S.~P. Yadav, ``Intuitionistic fuzzy data envelopment analysis: An
  application to the banking sector in india,'' {\em Expert Systems with
  Applications}, vol.~42, no.~11, pp.~4982--4998, 2015.

\bibitem{awadhmulti2022}
A.~Singh, S.~Yadav, and S.~Singh, ``A multi-objective optimization approach for
  dea models in a fuzzy environment,'' {\em Soft Computing},
  pp.~(https://doi.org/10.1007/s00500--021--06627--y), 2022.

\bibitem{wang2014fuzzy}
W.-K. Wang, W.-M. Lu, and P.-Y. Liu, ``A fuzzy multi-objective two-stage dea
  model for evaluating the performance of us bank holding companies,'' {\em
  Expert Systems with Applications}, vol.~41, no.~9, pp.~4290--4297, 2014.

\bibitem{chen2020comprehensive}
W.~Chen, S.-S. Li, J.~Zhang, and M.~K. Mehlawat, ``A comprehensive model for
  fuzzy multi-objective portfolio selection based on dea cross-efficiency
  model,'' {\em Soft Computing}, vol.~24, no.~4, pp.~2515--2526, 2020.

\bibitem{boubaker2020role}
S.~Boubaker, D.~T. Do, H.~Hammami, and K.~C. Ly, ``The role of bank affiliation
  in bank efficiency: a fuzzy multi-objective data envelopment analysis
  approach,'' {\em Annals of Operations Research}, pp.~1--29, 2020.

\bibitem{zamani2020position}
P.~Zamani, ``The position of multiobjective programming methods in fuzzy data
  envelopment analysis,'' {\em International Journal of Mathematical Modelling
  \& Computations}, vol.~10, no.~2 (SPRING), pp.~95--101, 2020.

\bibitem{wang2007measuring}
Y.-M. Wang, K.-S. Chin, and J.-B. Yang, ``Measuring the performances of
  decision-making units using geometric average efficiency,'' {\em Journal of
  the Operational Research Society}, vol.~58, no.~7, pp.~929--937, 2007.

\bibitem{zimmermann2011fuzzy}
H.-J. Zimmermann, {\em Fuzzy set theory—and its applications}.
\newblock Springer Science \& Business Media, 2011.

\bibitem{charnes1997data}
A.~Charnes, W.~Cooper, A.~Y. Lewin, and L.~M. Seiford, ``Data envelopment
  analysis theory, methodology and applications,'' {\em Journal of the
  Operational Research society}, vol.~48, no.~3, pp.~332--333, 1997.

\bibitem{marler2010weighted}
R.~T. Marler and J.~S. Arora, ``The weighted sum method for multi-objective
  optimization: new insights,'' {\em Structural and multidisciplinary
  optimization}, vol.~41, no.~6, pp.~853--862, 2010.

\bibitem{guo2001fuzzy}
P.~Guo and H.~Tanaka, ``Fuzzy dea: a perceptual evaluation method,'' {\em Fuzzy
  sets and systems}, vol.~119, no.~1, pp.~149--160, 2001.

\bibitem{wang2006centroids}
Y.-M. Wang, J.-B. Yang, D.-L. Xu, and K.-S. Chin, ``On the centroids of fuzzy
  numbers,'' {\em Fuzzy sets and systems}, vol.~157, no.~7, pp.~919--926, 2006.

\bibitem{zar2005spearman}
J.~H. Zar, ``Spearman rank correlation,'' {\em Encyclopedia of biostatistics},
  vol.~7, 2005.

\end{thebibliography}
	\newpage

\end{document}